\documentclass{amsart}

\newtheorem{thm}{Theorem}[section]
\newtheorem{cor}[thm]{Corollary}
\newtheorem{lem}[thm]{Lemma}
\newtheorem{conj}{Conjecture}

\theoremstyle{definition}
\newtheorem{defn}[thm]{Definition}
\theoremstyle{remark}

\numberwithin{equation}{section}

\def \eqskip { \vspace*{3mm} }

\newcommand{\RM}{\mathbb{R}}

\newcommand{\Aa}{\mathcal{A}^{\infty}}
\newcommand{\ag}{a(G)}

\newcommand{\di}{{\rm diag}}
\newcommand{\disps}{\displaystyle}
\newcommand{\ds}{\displaystyle}

\begin{document}

\title[Algebraic connectivity of graphs on surfaces]
{A Heawood--type result for the algebraic connectivity of graphs on surfaces}
\author{Pedro Freitas}
\thanks{Partially supported by FCT, Portugal}
%
%
\address{Departamento de Matem\'{a}tica,
Instituto Superior T\'{e}cnico, Av.Rovisco Pais, 1049-001 Lisboa,
Portugal.}
\email{pfreitas@math.ist.utl.pt}
\date{\today}

\begin{abstract}
We prove that the algebraic connectivity $a(G)$ of a graph embedded on a 
nonplanar surface satisfies a Heawood--type result. More precisely, it is shown 
that the algebraic connectivity of a surface $S$, defined as the supremum 
of $a(G)$ over all graphs that can be embedded 
in $S$, is equal to the chromatic number of $S$. Furthermore, and with the 
possible exception of the Klein bottle, we prove that this bound is 
attained only in the case of the maximal complete graph that can be 
embedded in $S$. In the planar case, we show that, at least 
for some classes of graphs which include the set of regular graphs, $a(G)$
is less than or equal to four. As an application of these results and 
techniques, we obtain a lower bound for the genus of Ramanujan graphs.

We also present some bounds for the asymptotic behaviour of $a(G)$ 
for certain classes of graphs as 
the number of vertices goes to infinity.
%
\end{abstract}

\maketitle
\section{Introduction}\label{intro}

In recent years, there has been some work relating the spectral radius 
of the adjacency matrix of a graph, $r(G)$, to its genus --
see~\cite{elzh,hong1,hong2}. The idea behind these results is to combine an
appropriate estimate for $r(G)$ in 
terms of the number of vertices and edges of the graph, with a direct consequence 
of Euler's formula giving an upper bound for the number of edges of a graph
embedded in a surface $S$. This allows the derivation of a bound for $r(G)$ in
terms on the number of vertices of $G$ and the genus of $S$.

In the continuous case there is also a relation between the
eigenvalues of a certain differential operator and the genus of that surface. More 
precisely, there exist upper bounds for the first nontrivial eigenvalue
of the Laplace-Beltrami operator on a surface in terms of its genus. 
The first of these results was obtained by Hersch for the case of the
sphere~\cite{hers}, and this was later generalized by Yang and
Yau to orientable surfaces of genus $\gamma$:
\begin{thm}\cite{yaya}
    Let $S$ be an orientable surface of genus $\gamma$ and let 
    $\lambda$ denote the first nontrivial eigenvalue of the 
    Laplace--Beltrami operator on $S$. Then
    \[
    \lambda \leq\frac{\disps 8\pi(1+\gamma)}{\disps A(S)},
    \]
    where $A(S)$ denotes the area of $S$.
\end{thm}
Note that in the case of general manifolds of dimension greater than 
two, it is 
known that no such results are possible~\cite{urak}.

These results suggest the derivation of similar bounds for the 
first nontrivial eigenvalue of the Laplacian operator defined on graphs. 
This quantity is related to the connectivity of a graph $G$, and was 
thus christened the algebraic connectivity of $G$, $a(G)$,
by Fiedler~\cite{fied}.

The main purpose of this paper is to study the maximum possible 
value of $a(G)$ for a graph which can be embedded on a surface of 
genus $\gamma$. It turns out that this value is equal to the Heawood 
number which appears in graph colouring problems~\cite{grtu,riyo,whbe}. Our main
result is that
for surfaces of positive genus this maximum possible value is given by the
algebraic connectivity of the 
maximal complete graph that is possible to embed on such a surface 
and that, with the possible exception of the Klein bottle, this value 
is attained only for this graph -- see
Section~\ref{secmain} for details. An immediate consequence is that the 
algebraic connectivity of a graph of genus $\gamma$ is bounded from 
above by the chromatic number of $S_{\gamma}$. In a sense, this is  
a surprising result since the algebraic connectivity and the chromatic 
number are not intimately related in general -- see Section~\ref{chrom}.

As in the case of map colouring problems, the techniques used to prove the
result in the case of positive genus do not apply to the case of the sphere, and
so this remains an open problem. We are, however, 
able to obtain some results under some restrictions which include the case of
planar regular graphs. These are
presented in Sections~\ref{planar} and~\ref{regsec}.

In Section~\ref{rama} we show how the results of the paper can 
be used to obtain estimates for the genus of some graphs. In 
particular, we apply this to the case of Ramanujan graphs.

These results show that the maximum value of the algebraic 
connectivity on a given surface is attained for a finite number of 
vertices. However, it is also of interest to know how
$a(G)$ behaves as the number of vertices becomes large.
A first step in this direction is done in Section~\ref{asymptsec} 
where we consider the supremum of the algebraic connectivity of certain
classes of graphs on a surface of fixed genus as
the number of vertices goes to infinity.

Finally, in Section~\ref{open} we
consider some open problems and 
present some conjectures.

\section{Preliminaries}\label{backg}

\subsection{Notation}

We begin by reviewing some concepts from graph theory that will be used in 
what follows. Let $G$ be a simple $n-$vertex graph ($n\geq3$), that is, a graph
with $n$ vertices and no loops nor multiple edges, and denote the sets of
vertices and edges by $V$ and $E$, respectively. We also denote the number of 
edges, $|E|$, by $e$. The degree $d_{i}$ of a vertex is the number of edges having that
vertex as one end, and the 
smallest and the largest of these numbers will be denoted by $d_{min}$ 
and $d_{max}$, respectively. 
The vertex connectivity $v(G)$ is defined as the minimal number of 
vertices whose removal (together with their adjacent edges) results in a 
disconnected graph. The girth of a graph is the length of the 
shortest cycle in the graph.

The adjacency matrix $A$ of a graph $G$ is defined by
\[
A=\left\{a_{ij}\right\}_{i,j=1}^{n}
\]
where $a_{ij}$ equals one if there is an edge connecting vertices $i$ and $j$,
and zero otherwise.

The Laplacian of a graph is defined to be the matrix $L=D-A$, where 
$D$ is the diagonal matrix $\di(d_{1},\ldots,d_{n})$. For some 
basic properties of this and related operators see, for 
instance,~\cite{chun,coli,fied}. The eigenvalues of the Laplacian will 
be denoted by
\[
0=\lambda_{1}\leq \lambda_{2}\leq\ldots\lambda_{n}.
\]
The second eigenvalue $\lambda_{2}$ is normally denoted by $a(G)$, and 
is called the algebraic connectivity of the graph.

As far as we are aware, the best bounds for $a(G)$ depending only on
the number of vertices and of edges of a graph remain those given by 
Fiedler in~\cite{fied}. The bound from that paper which will be of interest
here is contained in the following

\begin{thm}\cite{fied}\label{fiedthm}
    Let $G$ be a simple, connected, $n-$vertex graph with $e$ edges 
    and which is not complete. Then 
    \[
    \ag\leq v(G)\leq d_{min}\leq \frac{\ds 2e}{\ds n}.
    \]
\end{thm}

\subsection{Graphs on surfaces}

In what follows a surface $S$ is a compact connected $2-$manifold. From 
the classification of surfaces, we have that $S$ is either 
homeomorphic to a sphere with $h$ handles in the orientable case, or to 
the connected sum of $k$ projective planes in the non--orientable case. 
In this context, it is usual to denote the former by $S_{h}$ and the latter by $N_{k}$.
 The genus
$\gamma(S)$ of a surface $S$ is defined to be $1-\chi(S)/2$, where $\chi(S_{h})=2-2h$
and $\chi(N_{k})=2-k$ is the Euler characteristic of the surface.
Whenever a statement applies in 
both the orientable and non--orientable cases, we shall refer to the 
surface of genus $\gamma$ by $S_{\gamma}$, without stating explicitly 
which case is being considered.

An important aspect of topological graph theory
is the study of whether or not it is possible to embed a graph in a given 
surface. We say that there is an embedding of a graph $G$ in a surface $S$, if
there exist a one--to--one mapping of $V$ onto a set of $n$ distinct points in $S$, and a
mapping from $E$ to disjoint open arcs in $S$, such that no point in the image
of $V$ is contained in the image of an edge, and the image of an edge joining
two vertices is an arc joining the corresponding images -- see~\cite{whbe} for
this and other related concepts. If a graph $G$ is embedded in a surface $S$,
then the set $S\setminus G$ consists of a collection of connected components 
which are called the faces of $G$. If all the faces of an embedding
are homeomorphic to an open disk the graph is said to be cellularly 
embedded in $S$, and the embedding is called a $2-$cell embedding.
The orientable (non--orientable) genus of a graph is defined to be the smallest
possible genus of an orientable (resp. non--orientable) surface where $G$ is
embeddable. We shall refer to the orientable Euler characteristic of a
graph as the Euler characteristic of the orientable surface corresponding to the 
genus of the graph, and similarly for the non--orientable Euler 
characteristic. Where there is no danger of confusion, or when a 
result applies in both cases, we shall refer only to the genus or to the Euler
characteristic of a graph, without stating explicitly whether it 
refers to the orientable or non--orientable case.

The chromatic number of a surface $S$, 
$\kappa(S)$, is the maximum chromatic number of all graphs 
which can be embedded in $S$. In an analogous way, we define the 
algebraic connectivity of a surface as
\[
\mathcal{A}(S) = \sup_{G\in \mathcal{G}(S)}a(G),
\]
where $\mathcal{G}(S)$ is the set of all graphs that can be embedded 
in $S$.

Another question of interest 
is the study of the asymptotic behaviour of $a(G)$ for 
families of graphs embedded in a surface as the number of vertices goes to
infinity. Given an infinite family $\mathcal{F}$ of graphs embeddable on a
surface $S$ we define the asymptotic algebraic connectivity of $S$ by
\[
\Aa_{\mathcal{F}}(S) := \sup\left[\limsup_{n\to\infty} 
a(G_{n})\right],
\]
where $G_{n}\in\mathcal{F}$, $|G_{n}| = n$, and the supremum is taken 
over all possible such sequences -- when 
$\mathcal{F}=\mathcal{G}(S)$, we omit the subscript and write $\Aa(S)$. 

\subsection{Auxiliary results}
The basic result that allows us to relate the algebraic 
connectivity of a graph embedded in a surface to its genus is the following
consequence of Euler's formula:

\begin{thm}\label{eform}
    Let $G$ be a simple, connected, $n-$vertex graph with $e$ 
    edges and girth $g (\geq3)$. If $G$ is embeddable in a surface 
    $S$ with characteristic $\chi$, then
    \[
    e\leq \frac{\disps g}{\disps g-2}\left[n-\chi\right].
    \]
\end{thm}
For a proof see, for instance,~\cite{whbe}.

Combining this with Theorem~\ref{fiedthm} yields the following 
result, which will be used in the sequel.
\begin{cor}\label{asympt}
    Let $G$ be a non--complete graph of Euler characteristic $\chi$. Then
    \[
    a(G)\leq \frac{\ds 2g}{\ds g-2}\frac{\ds n-\chi}{\ds n}.
    \]
\end{cor}

We will also need a relation between the vertex connectivity of a 
graph and its genus. This is given by the following theorem, due to 
Cook.

\begin{thm}\cite{cook}\label{cookt}
    Let $G$ be a graph of nonpositive Euler characteristic $\chi$. Then
    \[
    v(G) \leq C(S_{\gamma}) := \left\lfloor\frac{\ds 
    5+\sqrt{49-24\chi}}{\ds 
    2}\right\rfloor.
    \]
\end{thm}
As mentioned in~\cite{plzh}, this still holds in the case of
the projective plane (non--orientable genus one). Cook has improved this result
when $G$ does not contain triangles. Here we shall only make use 
of this in the planar case.
\begin{thm}\cite{cook}\label{cookt2} Let $G$ be a planar graph of 
girth $g$. Then
\[
v(G)\leq \left\{
\begin{array}{ll}
5 & \mbox{ if } g=3,\\
3 & \mbox{ if } g=4,5,\\
2 & \mbox{ if } g\geq 6.
\end{array}
\right.
\]
\end{thm}

Finally, let $S$ be a surface of genus $\gamma$. Then the maximal 
complete graph of genus $\gamma$, which we denote by $K^{\gamma}$, is the 
complete graph on the largest possible number of vertices that can be embedded
in $S$.

A theorem of Ringel and Youngs gives the orientable genus of complete graphs.
\begin{thm}\cite{riyo}\label{maxgraphot}
    The orientable genus of the complete graph $K_{p}$ ($p\geq3$) is given by
    \[
    \gamma(K_{p}) = \left\lceil\frac{\ds (p-3)(p-4)}{\ds 
    12}\right\rceil.
    \]
\end{thm}
In the case of non--orientable genus, the corresponding result is due 
to Ringel.
\begin{thm}\cite{ring}\label{maxgraphnot}
    The non--orientable genus of the complete graph $K_{p}$ ($p\geq3$) is given by
    \[
    \gamma(K_{p}) = \left\lceil\frac{\ds (p-3)(p-4)}{\ds 
    6}\right\rceil,
    \]
    with the exception of $K_{7}$ for which we have $\gamma(K_{7}) = 3$.
\end{thm}


\section{A Heawood--type result for the algebraic connectivity\label{secmain}}

The main theorem of the paper is the following
\begin{thm}\label{maint}
    Let $G$ be a graph of genus $\gamma$ and nonpositive Euler 
    characteristic $\chi$. Then $a(G)\leq a(K^{\gamma})$.
    More precisely, and with the exception of the Klein bottle (non--orientable
    genus two), we have that
    \[
    \ag \leq a(K^{\gamma}) = H(S) := \left\lfloor\frac{\ds 
    7+\sqrt{49-24\chi}}{\ds 
    2}\right\rfloor,
    \]
    with equality if and only if $G=K^{\gamma}$

    In the case of the Klein bottle, 
    we have that $a(G)\leq a(K^{\gamma})= a(K_{6})=6<H(S)=7$.
    
    For planar graphs, $a(G)\leq 5$.
\end{thm}
The number $H(S)$ is Heawood's number for the surface $S$, and 
is related to the chromatic number of a surface -- see~\cite{grtu}, for 
instance.
As was pointed out in the Introduction, there is a great similarity between
this result and the corresponding result for the colouring of graphs 
known as the Heawood map--colouring problem. In fact,
a straightforward corollary to this theorem is that the algebraic
connectivity of a graph of genus $\gamma$ is bounded from above by the
chromatic number $\kappa(S_{\gamma})$.
\begin{cor}\label{cornonp}
    Let $G$ be a nonplanar graph of genus $\gamma$. Then 
    $a(G)\leq\kappa(S_{\gamma})$. Furthermore, if  $G\neq K^{\gamma}$ then
    $a(G)\leq a(K^{\gamma})-1$, except possibly in the case of the 
    Klein bottle.
\end{cor}
In terms of the algebraic connectivity of a surface, this may be 
stated as follows.
\begin{cor}
    For any surface $S$ of positive genus $\gamma$ we have that
    \[
    \mathcal{A}(S) = \kappa(S).
    \]
    Furthermore,  if  $G\neq K^{\gamma}$ is a graph embedded on a 
    surface $S$, then $a(G)\leq \mathcal{A}(S)-1$, except possibly in the case
    of the Klein bottle.
\end{cor}
 
\begin{proof}[Proof of Theorem~\ref{maint}]


Combining Theorem~\ref{cookt} with Fiedler's bound we have that for a 
nonplanar orientable graph $G$
\[
a(G) \leq v(G) \leq C(S) = H(S)-1.
\]
On the other hand, we have from Theorem~\ref{maxgraphot}
that for each value of $\gamma$ the complete graph 
$K_{H(S)}$, can be embedded in $S_{\gamma}$. Since 
$K_{H(S)}=H(S)$, we have that for noncomplete graphs
\[
a(G)\leq C(S)=H(S)-1=a(K_{H(S)})-1,
\]
which proves the result, as well as the second part of 
Corollary~\ref{cornonp}.

In the planar case the maximum possible
vertex connectivity is five, and thus $a(G)\leq 5$.

In the non--orientable case, and for $\gamma$ larger than two, we proceed in a
similar way to obtain the result for noncomplete graphs, now using 
Theorem~\ref{maxgraphnot}.

For the projective plane ($\gamma=1$), the result follows in the same 
way by using the fact that Theorem~\ref{cookt} extends to this case, and that
the maximal complete graph is now $K_{6}$.

In the case of the Klein bottle ($\gamma=2$), we have that the maximal 
complete graph that can be embedded there is $K_{6}$. From 
Theorem~\ref{cookt} we obtain that for noncomplete graphs $a(G)\leq 
v(G)\leq C(2) = 6 = a(K_{6})$, proving the result in 
this case. \end{proof}

\section{Planar graphs\label{planar}}

We will begin by obtaining some bounds on the algebraic connectivity which 
apply to general graphs. The main argument which will be used is a bound
based on a test 
function similar to that used in the proof of Cheeger's inequality -- 
see~\cite{coli}, for instance. The idea is that it 
should be possible to improve on bounds based on $v(G)$ in cases where 
there are clusters of vertices for which a sufficient large number of the edges
join vertices within the cluster. In order to do this, we need the following
definition.
\begin{defn} Let $G$ be a graph and $H$ be a proper nonempty subset of 
$V(G)$, formed by the vertices $y_{i}$, $i=1,\ldots,m$. The degree of the subset
$H$, is defined as
\[
d(H) = \sum_{i}^{m} \left(d_{i}-\tilde{d_{i}}\right),
\]
where $\tilde{d_{i}}$ is the number of edges joining two vertices 
$y_{i}$ and $y_{j}$, $1\leq i,j\leq m$.
\end{defn}

\begin{lem}\label{degreel}
    Let $G$ be a graph on $n$ vertices and let $H$ be a proper nonempty subset 
    of $V(G)$ with $m$ vertices and degree $d(H)$. Then
    \[
    a(G) \leq \frac{\ds d(H) n}{\ds m(n-m)}.
    \]
\end{lem}
\begin{proof}
    By the variational formulation for the eigenvalues of $G$, we have 
    that
    \begin{equation}\label{varform}
    a(G) \leq \frac{\ds \sum_{x\sim y} \left[ f(x)-f(y)\right]^{2}}
    {\ds \sum_{x\in G} f^{2}(x)},
    \end{equation}
    for all $f$ such that
    \[
    \sum_{i=1}^{n}f(x_{i}) = 0,
    \]
    and where $x\sim y$ denotes that the vertices $x$ and $y$ are adjacent.
    Let $f:V\to\RM$ be defined by
    \[
    f(x) = \left\{
    \begin{array}{ll}
	\frac{\ds n}{\ds m}-1, & x\in H\\
	-1, & x\in G\setminus H.
    \end{array}\right.
    \]
    We then have that $f$ is orthogonal to the vector with all entries 
    equal to $1$, and thus, by~(\ref{varform}), it follows 
    that
    \[
    \begin{array}{lll}
    a(G)&\leq & \frac{\ds d(H)\left(\frac{\ds n}{\ds m}\right)^{2}} {\ds 
    m\left(\frac{\ds n}{\ds m}-1\right)^2+n-m}\eqskip\\
    & = & \frac{\ds d(H)n^2}{\ds m^2\left( \frac{\ds n^2}{\ds 
    m}-n\right)}\eqskip\\
    &=& \frac{\ds d(H)n}{\ds m(n-m)}.
    \end{array}
    \]
\end{proof}

\begin{thm}\label{planarsmall} Let $G$ be a planar graph with $d_{max}$ smaller
    than or equal to five. Then $a(G)\leq 4$.
\end{thm}
\begin{proof} If $G$ does not have a triangle, then by 
Theorem~\ref{cookt2} $a(G)\leq 3$. Assume thus that $G$ has a 
triangle $T$. Then 
$d(T)\leq 3\times 5-6=9$ and hence, from Lemma~\ref{degreel}, it 
follows that
\[
a(G)\leq \frac{\ds 3n}{\ds n-3}.
\]
We thus have that $a(G)>4$ implies $n<12$. On the other hand, for a 
planar graph to have $a(G)$ greater than four, it must have $d_{min}$ 
equal to five. Since $d_{max}$ is less than or equal to five, from
Theorem~\ref{eform},
\[
\sum_{i=1}^{n} d_{i} = 5n = 2e\leq 6(n-2),
\]
and we obtain that for this to happen the graph
must have at least twelve vertices.
\end{proof}

\section{Regular graphs\label{regsec}}

Since for planar graphs $d_{min}$ is less than or equal to five,
Theorem~\ref{planarsmall} implies that for planar regular graphs $a(G)$
is less than or equal to four. We shall now consider the restriction to 
some particular cases of regular graphs.

We begin by obtaining a general bound for regular graphs (not 
necessarilly planar) of a given girth.

\begin{thm}\label{regasympt} Let $G$ be a regular graph of girth $g$ 
smaller than $n$ and Euler characteristic $\chi$. Then
\[
a(G) \leq \frac{\ds 2n}{\ds (n-g)(g-2)}\left( 2-\frac{\ds g\chi}{\ds n}
\right).
\]
\end{thm}
\begin{proof} 
Take $H$ in Lemma~\ref{degreel} to correspond to the cycle $C_{g}$. Then $d(H)
= g(d-2)$ and
\[
\begin{array}{lll}
a(G)& \leq & \frac{\ds (d-2)n}{\ds n-g}\eqskip\\
& = & 2\left( \frac{\ds e}{\ds n}-1\right)\frac{\ds n}{\ds n-g}\eqskip\\
& \leq & 2\left( \frac{\ds g}{\ds g-2} \frac{\ds (n-\chi)}{\ds n}-1\right)
\frac{\ds n}{\ds n-g}\eqskip\\
& = & \frac{\ds 2}{\ds g-2}\left( 2-\frac{\ds g\chi}{\ds n}\right)\frac{\ds n}{\ds 
n-g}.
\end{array}
\]
\end{proof}
Note that for a fixed girth and sufficiently large $n$ this is better than
Corollary~\ref{asympt} which is obtained from Fiedler's bound.

In the case of planar graphs $\chi$ equals two and this bound takes the following
simple form.

\begin{cor} \label{plgi} Let $G$ be a regular planar graph of girth $g$ smaller 
than $n$. Then
\[
a(G)\leq \frac{\ds 4}{\ds g-2}.
\]
\end{cor}

When the girth is equal to $n$, we have the case of $2-$regular 
graphs (the cycle $C_{g}$) for which it is easy to see that, with the 
exception of $K_{3}$, the algebraic connectivity is always less than 
or equal to two, this value being attained for the cycle $C_{4}$. A 
similar result holds for planar cubic graphs.
\begin{thm} Every planar cubic graph other than $K_{4}$ has
    \[
    a(G)\leq 2.
    \]
\end{thm}
\begin{proof} If the girth $g$ is larger than or equal to four, the result 
follows from the bound in Corollary~\ref{plgi}. Assume thus that $g$ 
equals three. Then $G$ has a triangle and proceeding as in the proof of 
Theorem~\ref{degreel} but now using the fact that the degree is three, 
we obtain that
\[
a(G)\leq \frac{\ds n}{\ds n-3}.
\]
This will be greater than two provided that $n$ is less than $6$. 
Since we are interested in cubic graphs, the case of $n$ equal to five 
is excluded and the only cubic graph when $n$ equals four is $K_{4}$.
\end{proof}
The result is sharp, in the sense that the bound is attained by both 
the $3-$prism and the cube graphs.

\section{Algebraic connectivity and chromatic numbers\label{chrom}}

In general, and without any further assumptions, we cannot expect a deep relation between $a(G)$ and 
$\kappa(G)$, in the sense that there exist graphs for which $a(G)<\kappa(G)$ (any 
connected graph with $a(G)$ less than two), graphs for which $a(G)=\kappa(G)$ 
(complete graphs), and graphs for which $a(G)>\kappa(G)$ (any 
complete bipartite graph $K_{p,q}$ with $p$ and $q$ both greater than 
two). However, it is possible to prove the following

\begin{thm} \label{chromt} Let $G$ be a noncomplete graph on $n$ vertices and
    with chromatic number $\kappa(G)$. Then
\[
a(G) \leq n - \left\lceil \frac{\ds n}{\ds \kappa(G)}\right\rceil.
\]
\end{thm}
\begin{proof} The graph $G$ can be divided into subsets 
$X_{i}, i=1,\ldots,\kappa(G)$, such that there are no edges connecting 
vertices within each set $X_{i}$. This means that the complementary 
graph of $G$, $G^{c}$, contains the complete graphs $K_{n_{i}}$, 
where $n_{i}=|X_{i}|, i=1,\ldots\kappa(G)$. We also have that at 
least one of the numbers $n_{i}$ is greater than $\lceil 
n/\kappa(s)\rceil$ and, since $G$ is not 
complete, greater than one.  Hence 
\[\lambda_{n}(G^{c})\geq\left\lceil \frac{\ds n}{\ds 
\kappa(G)}\right\rceil,\] from which the result follows, since
$a(G)\leq n-\lambda_{n}(G^{c})$.
\end{proof}

In the special case of bichromatic graphs of a given Euler characteristic it is
possible to improve on this result.

\begin{thm}\label{bichromt} If $G$ is a graph on $n$ vertices with $\kappa(G)=2$
    and Euler characteristic $\chi$, then
    \[
    a(G) \leq 4 \frac{\ds n-\chi}{\ds n}.
    \]
\end{thm}
\begin{proof}
    Since $G$ is bichromatic, its girth must be greater than three. 
    The result now follows from Corollary~\ref{asympt}.
\end{proof}
As a consequence, we obtain that bichromatic planar graphs satisfy $a(G)<4$.
However, in this case it is possible to obtain a better bound using
Theorem~\ref{cookt2}
\begin{thm} If $G$ is a bichromatic planar graph, then $a(G)\leq 3$.
\end{thm}
\begin{proof} It follows as before but now using Theorem~\ref{cookt2} 
instead of Corollary~\ref{asympt}.
\end{proof}

It is, of course, also possible to use the results in~\cite{cook} for the case
of general orientable genus and where the graph's girth is greater than or equal
to four to obtain further bounds for bichromatic graphs. However, these will
not be as good as those in Theorem~\ref{bichromt}, at least for $\gamma$ larger
than one and sufficiently large $n$.

In general, bichromatic graphs can have an arbitrarilly high 
algebraic connectivity, as can be seen from the case of the complete 
bipartite graph $K_{p,q}$, for which $\kappa(K_{p,q})$ equals 
two, while $a(K_{p,q})=\min\{p,q\}$. Note that when $p$ equals $q$ this 
gives equality in the case of Theorem~\ref{chromt}.

On the other hand, it is also possible to keep $a(G)$ bounded while 
making the chromatic number as large as desired. To see this, consider 
the graph $G$ obtained from $K_{n}$ by adding a vertex which is 
connected to one single vertex in $K_{n}$, that is, 
$G=(K_{n-1}\cup \{x\})\dot{+}K_{1}$. The spectrum of this graph is 
given by $(0,1,n,\ldots,n,n+1)$ (and thus $a(G)=1$), while its
chromatic number grows with $n$.

These two examples show that unless we impose an extra restriction 
such as fixing the genus of the graph, we should not expect a close relation
between the algebraic connectivity and the chromatic number of a graph.

\section{\label{rama}A lower bound for the genus of Ramanujan graphs}

It is possible to use the bounds in Corollary~\ref{cornonp} 
to obtain estimates for the genus of a given graph, provided one has 
a lower bound for $a(G)$ -- this is actually just a consequence from 
the fact that $a(G)$ is less than or equal to $C(S)$. Here we apply this to the
case of Ramanujan graphs -- see~\cite{chun,lps}.

\begin{thm} Let $G$ be a Ramanujan graph on $n$ vertices and of 
degree $d$, with $9\leq d\neq n-1$. Then its orientable genus satisfies
\[
\gamma\geq \left\lceil\frac{\ds (2d-4\sqrt{d-1}-5)^{2}-1}{\ds 
48}\right\rceil.
\]
\end{thm}
\begin{proof} The algebraic connectivity of a Ramanujan graph 
satisfies (see~\cite{chun}, page $97$, for instance)
\[
a(G)\geq d-2\sqrt{d-1}.
\]
Since by Corollary~\ref{cornonp} the algebraic connectivity of a noncomplete
graph which can be embedded on a surface $S$ satisfies
\[
a(G)\leq a(K^{\gamma})-1=H(S)-1=C(S),
\]
it follows that
\[
\sqrt{49-24\chi}\geq 2d-4\sqrt{d-1}-5.
\]
From this we see that if $d$ is smaller than $9$ then nothing can be 
concluded, while for $d$ greater than or equal to $9$ we obtain the desired result.
\end{proof}
A similar result can be obtained for the case of the nonorientable 
genus in the same way, except that then one has to consider the 
exceptional case of genus two separately.

\section{Asymptotic behaviour of $a(G)$\label{asymptsec}}

From Corollary~\ref{cornonp} we have that the algebraic connectivity of a 
graph embedded in a surface $S$ is bounded from above independently 
of the number of vertices of the graph and, in fact, if $S$ is neither the
sphere nor the Klein bottle, then $\Aa(S)\leq \mathcal{A}(S)-1$. 
A better bound for $\Aa(S)$ is given by Corollary~\ref{asympt} which yields that
$\Aa(S)$ is less than or equal to $6$, independently of $S$.

On the other hand, by considering the sequence of (planar) {\it double wheel}
graphs, that is, $G_{n}=C_{n}\dot{+}2K_{1}$, for which
\[
a(G_{n}) = 
\min\left\{a(C_{n})+2,a(2K_{1})+n\right\} = 
\min\left\{4-2\cos(2\pi/n),n\right\},
\]
we obtain that $\Aa(S)$ is 
greater than or equal to two and so we have that
\[
2\leq \Aa(S)\leq 6.
\]

We shall now restrict our attention to certain classes of graphs. An
immediate consequence of Theorem~\ref{regasympt} is a bound for the asymptotic 
connectivity of regular graphs with a given fixed girth $g$.

\begin{thm}\label{corasympt1}
    Let $\mathcal{G}_{r,g}(S)$ be the set of all regular graphs with 
    girth equal to a fixed number $g$ that can 
    be embedded in a given surface $S$. Then
    \[
    \Aa_{\mathcal{G}_{r,g}}(S)\leq \frac{\ds 4}{\ds g-2}.
    \]
\end{thm}
Note that the dependence of the bound in Theorem~\ref{regasympt}
is not monotonic in $g$ and so we cannot obtain directly a uniform bound which 
would allow us to conclude that the asymptotic connectivity of 
regular graphs is less than or equal to four. However, it is 
possible to prove this combining Theorem~\ref{regasympt} with 
Corollary~\ref{asympt}.
\begin{thm}\label{asymptreg}
    Let $\mathcal{G}_{r}(S)$ be the set of all regular graphs 
    that can be embedded in $S$. Then
    \[
    \Aa_{\mathcal{G}_{r}}(S)\leq 4.
    \]
\end{thm}
\begin{proof}
    If $g$ is equal to three, Theorem~\ref{regasympt} gives that 
    \[
    a(G)\leq \frac{\ds 2n}{\ds n-3}\left(2-\frac{\ds 3\chi}{\ds 
    n}\right)=\frac{\ds 4n-6\chi}{\ds n-3}.
    \]
    On the other hand, 
    if the girth is larger than three, then by Corollary~\ref{asympt} it 
    follows that
    \[
    a(G)\leq \frac{\ds 4(n-\chi)}{\ds n}.
    \]
    Since both bounds converge to four as $n$ goes to infinity, the 
    result follows.
\end{proof}
Combining this with a bound of Alon and Boppana, we 
obtain the following bound for the specific case of $d-$regular graphs 
for a fixed integer $d$.
\begin{cor}
    For a given integer $d$ greater than or equal to two, let
    $\mathcal{G}^{d}(S)$ denote the set of $d-$regular graphs that can be 
    embedded in $S$. We have that
    \[
    \Aa_{\mathcal{G}^{d}}(S)\leq
    \left\{
    \begin{array}{ll}
	d-2\sqrt{d-1}, & d\leq 10\\
	4, & d>10.
    \end{array}
    \right.
    \]
\end{cor}
\begin{proof}
    From a result by Alon and Boppana quoted in~\cite{alon} -- see 
    also~\cite{chun} --, we have 
    that for a $d-$regular graph
    \[
    a(G)\leq d-2\sqrt{d-1}+O(\log_{d}n)^{-1},
    \]
    from which it follows that $\Aa_{\mathcal{G}^{d}}(S)\leq 
    d-2\sqrt{d-1}$. This is smaller than four for $d$ smaller then or 
    equal to $10$. For $d$ larger than $10$, we use the bound from 
    Theorem~\ref{asymptreg}.
\end{proof}

\section{Discussion\label{open}}

As in the proof of Theorem~\ref{planarsmall}, we see that for $a(G)$ to be
greater than four then a planar graph must have at least twelve vertices, since this
is a condition for the minimum degree of a planar graph to be equal to 
five. We thus have that $a(G)\leq 4$ both when $n$ is smaller than 
twelve or when $d_{max}$ is smaller than or equal to five.
Both this and the similarity with colouring problems suggest that part of
Theorem~\ref{maint} extends to planar graphs, that is,
that the maximal algebraic connectivity of a planar graph is four, the algebraic 
connectivity of $K_{4}$, but we haven't been able to 
prove it. However, if that is the case, then the maximal algebraic 
connectivity will not be uniquely attained in this case, since the 
octahedron (the join of
$2K_{1}$ with the cycle on four vertices, $2K_{1}\dot{+}C_{4}$), has as its 
spectrum $(0,4,4,4,6,6)$. We have the following conjecture
\begin{conj} If $G$ is a planar 
graph, then $a(G)\leq 4$, with equality if and only if $G=K_{4}$ or 
$G=2K_{1}\dot{+}C_{4}$. Furthermore, if $G$ is neither of these 
graphs, then $a(G)$ is less than or equal to three.
\end{conj}

Regarding chromatic numbers, note that $\kappa(2K_{1}\dot{+}C_{4})$ is three,
and so it is possible to
have a planar graph whose chromatic number is three while its
algebraic connectivity equals four, although if the conjecture above 
holds, then this will be the only graph for which this will happen.

As we have seen, for bichromatic planar graphs we must have that 
$a(G)$ can be at most three. In fact, we conjecture that this value 
can be improved.

\begin{conj} Every planar bichromatic graph has $a(G)\leq2$.
\end{conj}

In the case of the Klein bottle, and since it is possible to have
equality between $C(N_{2})$ and $a(K_{6})$, we have not been able to prove that
equality holds only for the complete graph $K_{6}$. However, we 
believe this to be the case.

\begin{conj} For the Klein bottle (the non--orientable case of genus 
two), $a(G)$ is $6$ if and only if $G=K_{6}$.
\end{conj}

Finally, regarding the upper bound for the asymptotic algebraic connectivity 
in the general case, we remark that
this is most likely far from being optimal, since we have used 
methods which are based mainly on local properties of graphs. We conjecture
that the value of the asymptotic algebraic connectivity in the general case is
independent of the (fixed) surface $S$, and that it will in fact be equal 
to two, the algebraic connectivity of the double wheel graph.
\begin{conj} $\Aa(S) \equiv 2$.
\end{conj}

\end{document}